\documentclass[12pt]{amsart}
\usepackage[parfill]{parskip}    
\usepackage{amsfonts}
\usepackage{amsmath}
\usepackage{amsthm}
\usepackage{epstopdf}
\usepackage{hyperref}
\usepackage[all]{xy}
\usepackage[pdftex]{color,graphicx}
\usepackage{psfrag}
\usepackage{times}

\title{Examples of planar tight contact structures with support norm one}
\author{Tolga Etg\"u}
\address{Department of Mathematics, Ko\c{c} University, \.Istanbul 34450 TURKEY \newline Mathematical Sciences Research Institute, Berkeley, CA 94720 USA}
\email{tetgu@ku.edu.tr}
\author{Yank\i\ Lekili}
\address{Mathematical Sciences Research Institute, Berkeley, CA 94720 USA}
\email{ylekili@msri.org}
\thanks{We would like to thank Selman Akbulut and Burak Ozbagci for helpful conversations. TE was partially supported by the Scientific and Technological Research Council of Turkey.}

\newtheorem*{thm*}{Theorem} \theoremstyle{definition}

\newcommand{\QED}{\vspace{-.23in}\begin{flushright}\qed\end{flushright}}

\begin{document}

\begin{abstract}

We exhibit an infinite family of tight contact structures with the property that none
of the supporting open books minimizes the genus and maximizes the Euler characteristic
of the page simultaneously, answering a question of Baldwin and Etnyre in~\cite{EB}. 

\end{abstract} \maketitle

Let $Y$ be a closed oriented 3--manifold and $\xi$ be a contact structure on $Y$.
Recall that an open book is a fibration $\pi: Y - B \to S^1$ where $B$ is an oriented
link in $Y$ such that the fibres of $\pi$ are Seifert surfaces for $B$. The contact
structure $\xi$ is said to be supported by an open book $\pi$ if $\xi$ is the kernel
of a one-form $\alpha$ such that $\alpha$ evaluates positively on the positively
oriented tangent vectors of $B$ and $d\alpha$ restricts to a positive area form on each fibre
of $\pi$. It is well known that every contact structure $\xi$ is supported by an open book
on $Y$ and all open book decompositions of $Y$ supporting $\xi$ are equivalent up to
positive stabilizations and destabilizations \cite{G}, but given a contact 3--manifold $(Y, \xi)$ it is not 
always easy to find a ``simple" supporting open book. One natural measure of simplicity comes from the Euler characteristic of a page which is decreased by stabilization. The genus of a page is another useful indicator of simplicity. 

In \cite{EO}, Etnyre and Ozbagci define three numerical
invariants of $\xi$, called the \emph{support norm, support genus} and \emph{binding
number}, respectively, in terms of its supporting open books: \[ \text{sn}(\xi) =
\text{min}\{ -\chi(\pi^{-1}(\theta))| \pi : Y-B \to S^1 \text{\ supports\ } \xi \} \]
\[ \text{sg}(\xi) = \text{min}\{ g(\pi^{-1}(\theta))| \pi : Y-B \to S^1 \text{\
supports\ } \xi \} \] \[ \text{bn}(\xi) = \text{min}\{ |B| | \pi : Y-B \to S^1 \text{\
supports\ } \xi \text{\ and \ } g(\pi^{-1}(\theta)) = sg(\xi) \} \ , \] where $\theta$
is any point in $S^1$ , $g(.)$ is the genus, and $|.|$ is the number of components.
In general, these invariants are hard to compute. It is known that $\text{sg}(\xi) =
0$ if $\xi$ is overtwisted \cite{E}, and in general there are obstructions for a
contact structure to have support genus zero (\cite{E, OSS}). However, there is no
known example of a contact structure with support genus greater than one. Even if
$\xi$ is overtwisted, it is not easy to determine $\text{bn}(\xi)$. Furthermore, it is
known that no two of these invariants determine the third \cite{EB}.    

It is obvious from the above definitions that \[ \text{sn}(\xi) \leq 2 \text{sg}(\xi) +
\text{bn}(\xi)-2 \ ,  \]  and that equality holds when $\text{bn}(\xi) \leq 3$. 

In \cite{EB}, Baldwin and Etnyre exhibit examples of {\it overtwisted}
contact structures which make the above inequality strict and ask whether the inequality
can be strict for tight contact structures. Here we give an infinite family of tight
contact structures (exactly one of which is Stein fillable) for which this inequality
is strict. 

\begin{figure}[!h]
\centering 
\includegraphics[scale=0.6]{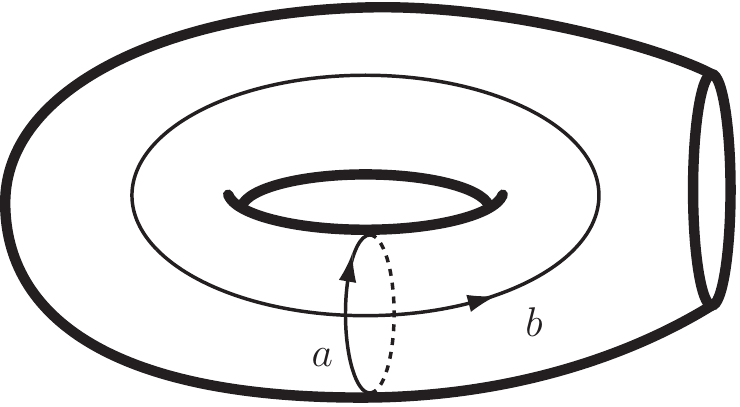}  \caption{}
\label{Figure0} \end{figure}

Let $T_0$ be genus one surface with one boundary component and consider the family of
diffeomorphisms $\phi_m = (\tau_a\tau_b)^3 \tau_a^{-m-4} $, for $m\geq 0$, where $a$
and $b$ are simple closed curves given in Figure \ref{Figure0} 
and $\tau$ denotes the right-handed Dehn twist along the
corresponding curve. For later use, we orient $a$ and $b$ so that $a\cdot b= -1$. Let
$(Y_m, \xi_m)$ denote the contact manifold supported by the open book decomposition
$(T_0, \phi_m)$. 

\begin{thm*}

The contact structure $\xi_0$ is Stein fillable and $\xi_m$ is tight but not
Stein fillable for $m>0$. Furthermore,
\[ \text{sn}(\xi_m) = 1 \ , \  \text{sg}(\xi_m)=0 \ , \ 3<\text{bn}(\xi_m) \leq m+5 \] 

In particular, $\text{sn}(\xi_m) < 2 \text{sg}(\xi_m) + \text{bn}(\xi_m)-2 $. 

\end{thm*} 
\emph{Proof.} The fact that every $\xi_m$ is tight with nontrivial Heegaard Floer invariant $c(Y_m, \xi_m) \in \widehat{HF}(-Y_m, \mathfrak{s}_{\xi_m})$ follows from Theorem $4.3$ in \cite{HKM2}.
Now, using the relations $\tau_a \tau_b \tau_a = \tau_b \tau_a \tau_b$ and $\tau_{f(\gamma)} = f \tau_{\gamma} f^{-1}$ for any automorphism $f$ and simple closed curve $\gamma$, 
\[ \phi_0 = (\tau_a\tau_b)^3 \tau_a^{-4} 
        = \tau_a^{-2} \tau_a\tau_b\tau_a\tau_b\tau_a\tau_b\tau_a^{-2} 
        = \tau_a^{-1}\tau_b\tau_a \tau_a\tau_b\tau_a^{-1}
        = \tau_{a+b}\tau_{a-b}  \]

\begin{figure}[!h]
        \centering
        \includegraphics[scale=0.5]{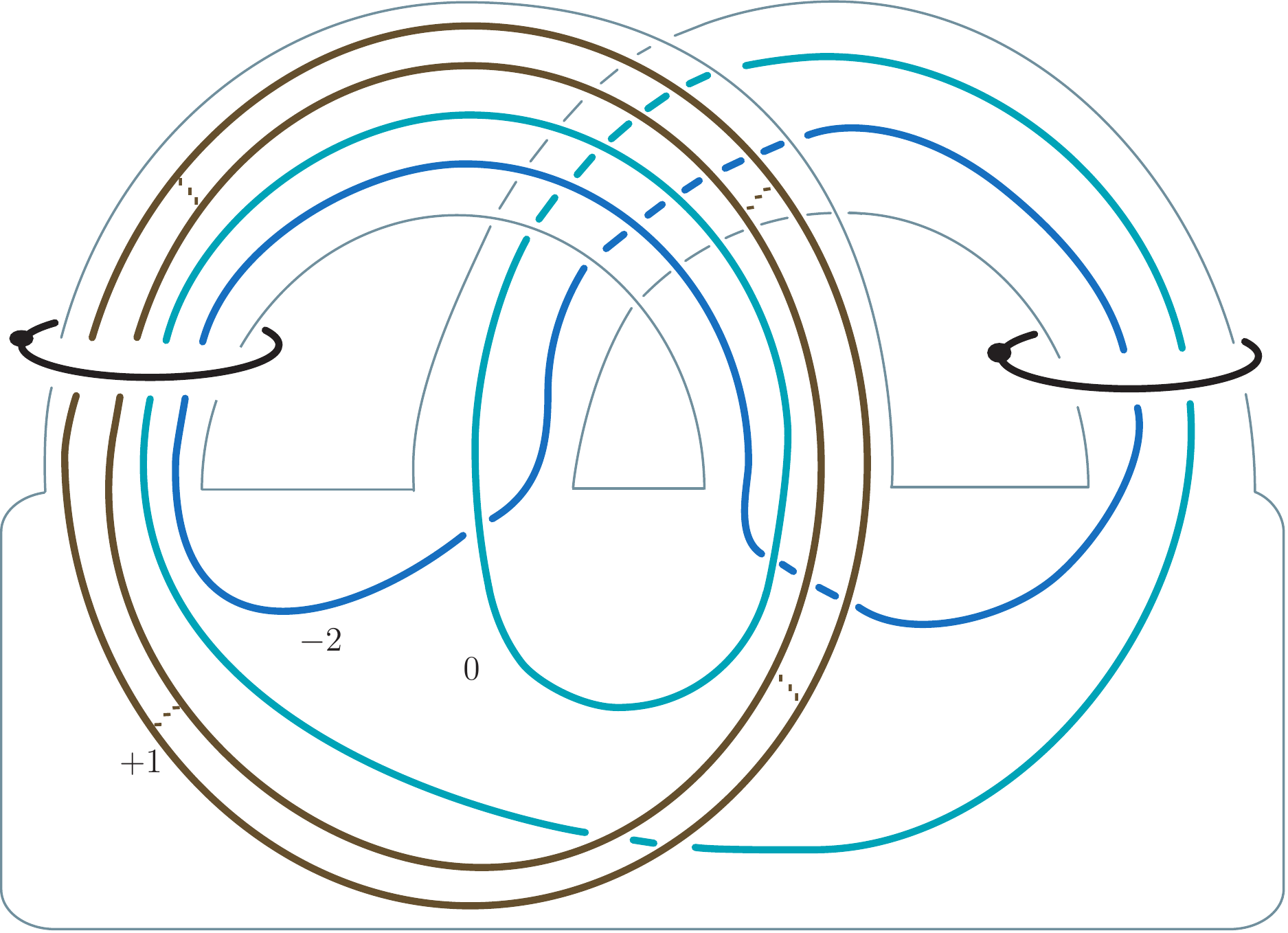}  \caption{Handlebody diagram of $X_m$ with two 1-handles, $m$ $
        +1$-framed 2-handles, a $-2$-framed 2-handle, and a $0$-framed 2-handle}
\label{Figure1}
\end{figure}

Since it is supported by an open book whose monodromy is a
product of right-handed Dehn twists, $\xi_0$ is Stein fillable. In general, we have $\phi_m  = \tau_{a+b}\tau_{a-b}\tau_{a}^{-m}$. Using this factorization, we draw a handlebody diagram of a 4--manifold
$X_m$ with boundary $Y_m$ in Figure \ref{Figure1}. 

\begin{figure}[!h]
\centering
\includegraphics[scale=0.9]{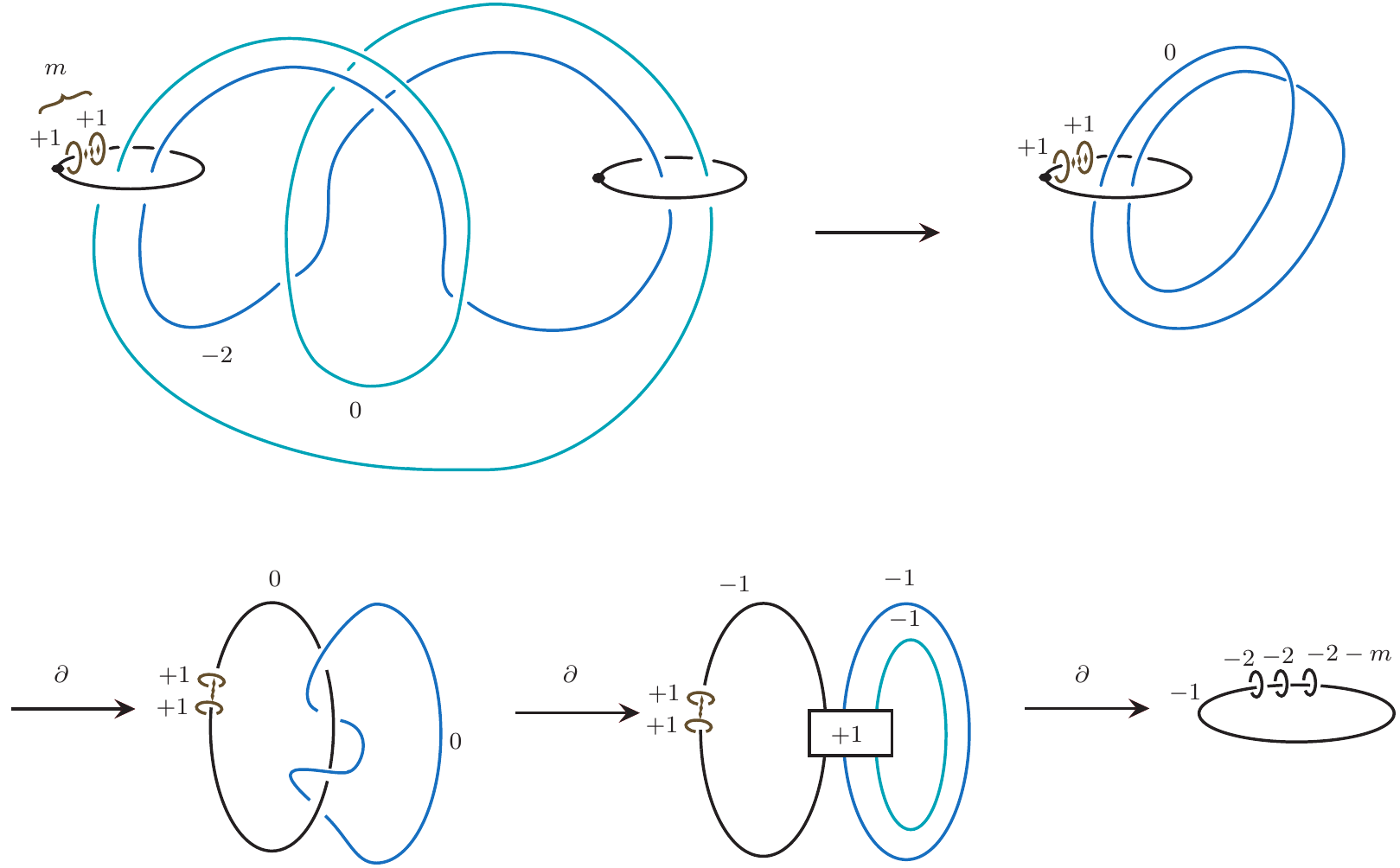} \caption{Seifert fibered 3--manifold description of $Y_m$}
\label{Figure2}
\end{figure}

\begin{figure}[!h]
\centering
\includegraphics[scale=1.1]{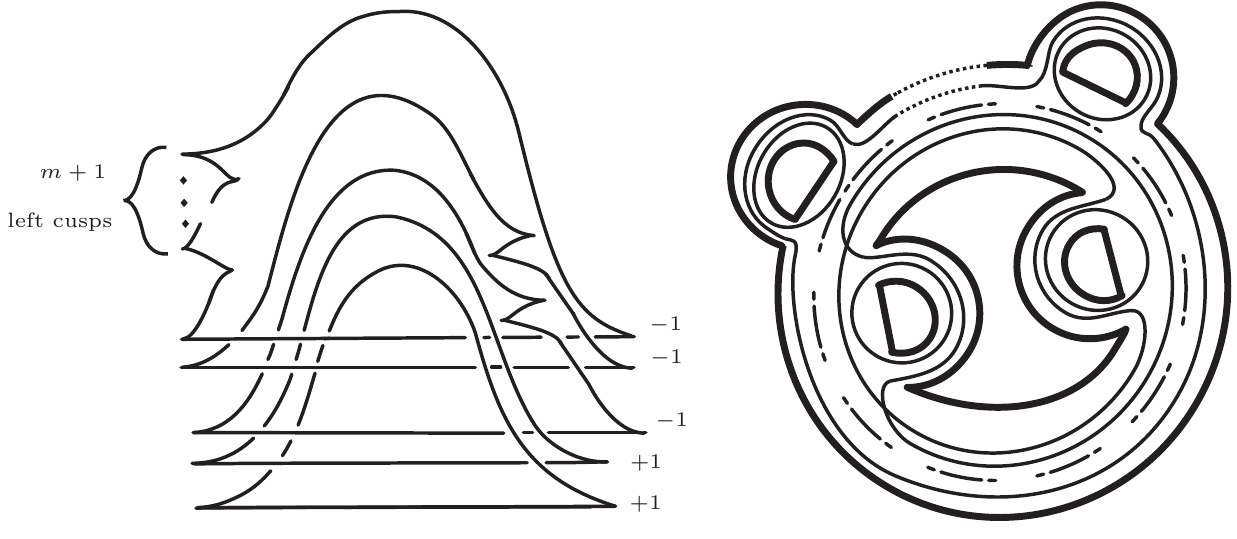} \caption{On the left: Contact surgery diagram of $\xi_m$ as given by Figure 7 in \cite{GLS}. On the right: A planar open book supporting $\xi_m$ with $m+5$ boundary components, the monodromy is negative Dehn twist around the middle dashed curve and positive Dehn twist around all the other curves.}
\label{Figure3}
\end{figure}

Figure \ref{Figure2} describes a way to see that $Y_m$ is diffeomorphic to the Seifert fibered 3--manifold
$M(-1;\frac{1}{2} , \frac{1}{2} , \frac{1}{m+2})$. 
A complete classification of tight contact structures on
$M(-1;\frac{1}{2},\frac{1}{2} , \frac{1}{m+2})$ is given in Section 4 of \cite{GLS}.
From their classification, it follows that all the tight contact structures on these manifolds are supported by
planar open books, i.e. $\text{sg}(\xi_m)=0$. In fact, we can pinpoint precisely
the contact isotopy class of $\xi_m$ from this classification by calculating a Hopf
invariant, $d_3(\xi_m)$. Indeed, in Theorem 4.1 and Proposition 5.1 of \cite{B} Baldwin shows that $Y_m$ is an L-space and calculates the correction term $d(Y_m,
\mathfrak{s}_{\xi_m}) = -\frac{m}{4}$. Since we also know that $c(Y_m,\xi_m)$ is
non-zero, and it has grading equal to $-d_3(\xi_m)$ \footnote{Here we follow the convention in \cite{GLS} where the Hopf invariant is shifted by $1/2$ so that it is $0$ for the standard contact structure on $S^3$. } by Proposition 4.6 of \cite{OS}, we conclude that
$d_3(\xi_m)=d(Y_m,\mathfrak{s}_{\xi_m})=-\frac{m}{4}$. (Note that this calculation
can also be done by drawing a contact surgery diagram associated with $\phi_m$).
For each $m>0$ there are three tight contact structures on $Y_m$ \cite{GLS} exactly one of which has $d_3$ invariant equal to $-\frac{m}{4}$, and it is given by the contact surgery diagram in Figure~\ref{Figure3}.
Note that the fact that $\xi_m$ is Stein fillable if and only if $m=0$ follows from
Theorem 4.13 in \cite{GLS}. 

\begin{figure}[!h]
\centering
\includegraphics[scale=0.6]{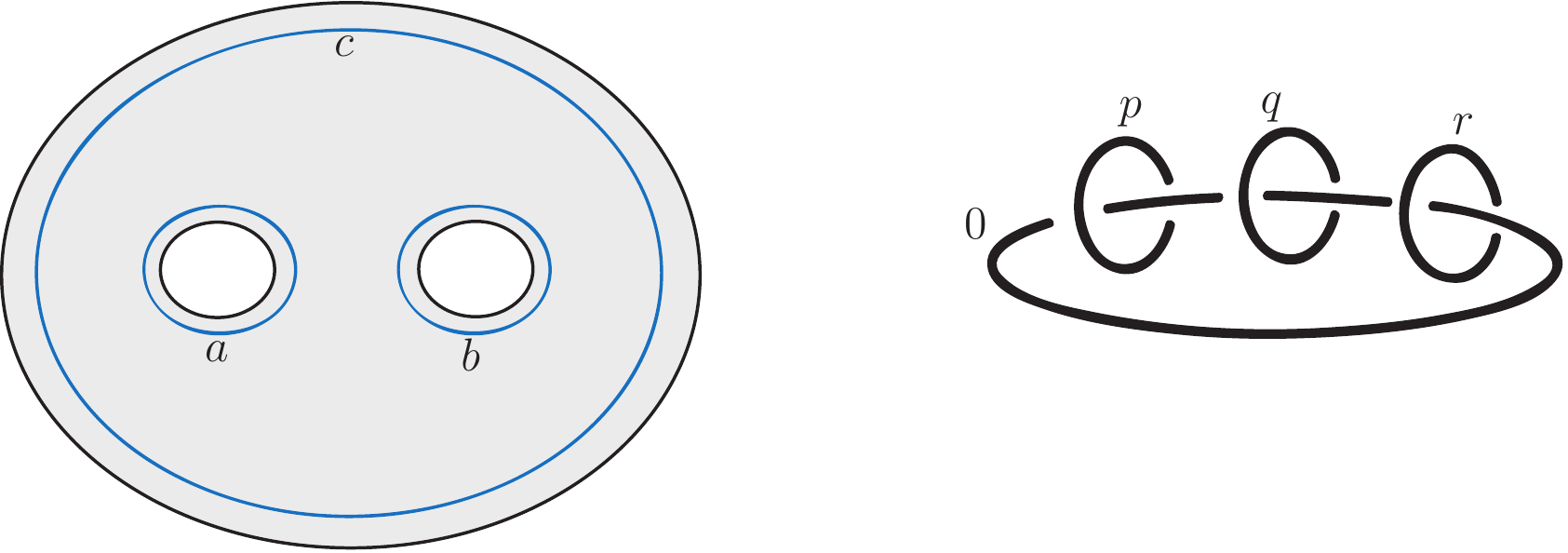}   \caption{A surgery picture of the 3--manifold given by the planar open book with three binding components and monodromy $\phi = \tau_a^p \tau_b^q \tau_c^r$, where $\tau$ denotes the right-handed Dehn twist along the corresponding curve}
\label{Figure4}
\end{figure}

So far, we have shown that $\text{sg}(\xi_m)=0$, and $\text{sn}(\xi_m)\leq 1$, where
the latter follows because we started with an open book supporting $\xi_m$
with pages a genus one surface with one boundary component. Furthermore, Figure
\ref{Figure3} gives a planar open book supporting $\xi_m$ with $m+5$ boundary components, hence $\text{bn}(\xi_m)\leq m+5$. Next, observe that
$\text{sn}(\xi) < 1 $ implies that $\xi$ is a contact structure on a lens space. Hence
$\text{sn}(\xi_m) =1$. To finish, we need to show that $\text{bn}(\xi_m) \neq 3$. Any 3--manifold with a planar open book with three binding components is given by
a surgery diagram as in Figure \ref{Figure4}. 
These are connected sums of lens spaces if $\{ 0,\pm 1 \} \cap \{p,q,r \} \neq \emptyset$, and small Seifert fibered spaces with $e_0=\lfloor -\frac{1}{p} \rfloor + \lfloor -\frac{1}{q} \rfloor + \lfloor -\frac{1}{r} \rfloor$ otherwise. Since $e_0 (Y_m)= -1$ any open book decomposition of $Y_m$ with planar pages and three binding components must have exactly two of $p,q$ and $r$ negative, and in that case the monodromy is not
right-veering. Therefore, these open books cannot support the tight contact
structures $\xi_m$  by \cite{HKM}. \QED

\end{document}